\documentclass[11pt]{article}
\usepackage{amsmath}
\usepackage{dsfont}
\usepackage{mathrsfs}
\usepackage{amsmath,amssymb}
\usepackage{amsfonts}
\usepackage{hyperref}
\usepackage{amsthm}
\usepackage{graphicx}
\usepackage{subfigure}
\usepackage{xcolor}

\hfuzz=\maxdimen
\theoremstyle{definition}
\newtheorem{lemma}{Lemma}[section]

\newtheorem{theorem}[lemma]{Theorem}
\newtheorem{corollary}[lemma]{Corollary}

\newtheorem{remark}{Remark}

\numberwithin{equation}{section}

\DeclareFixedFont{\Acknowledgment}{OT1}{cmr}{bx}{n}{14pt}
\begin{document}

\title{\bf A $p$-th Yamabe equation on graph}
\author{Huabin Ge}
\maketitle

\begin{abstract}
Assume $\alpha\geq p>1$. Consider the following $p$-th Yamabe equation on a connected finite graph $G$:
$$\Delta_p\varphi+h\varphi^{p-1}=\lambda f\varphi^{\alpha-1},$$
where $\Delta_p$ is the discrete $p$-Laplacian, $h$ and $f>0$ are fixed real functions defined on all vertices. We show that the above equation always has a positive solution $\varphi$ for some constant
$\lambda\in\mathds{R}$.
\end{abstract}


\section{Introduction}
The well known smooth Yamabe problem asks for the considering of the following smooth Yamabe equation \cite{Aubin,Lee,Yamabe}
\begin{equation*}
\Delta\varphi+h(x)\varphi=\lambda f(x)\varphi^{N-1}
\end{equation*}
on a $C^{\infty}$ compact Riemannian manifold $M$ of dimension $n\geq 3$, where $h(x)$ and $f(x)$ are $C^{\infty}$ functions on $M$, with $f(x)$ everywhere strictly positive and $N=2n/(n-2)$. The problem is to prove the existence of a real number $\lambda$ and of a $C^{\infty}$ function $\varphi$, everywhere strictly positive, satisfying the above Yamabe equation. In this short paper, we consider the corresponding discrete Yamabe equation
\begin{equation*}
\Delta\varphi+h\varphi=\lambda \varphi^{\alpha-1},\;\;\alpha\geq2
\end{equation*}
on a finite graph. More generally, we shall establish the existence results of the following $p$-th discrete Yamabe equation
\begin{equation*}
\Delta_p\varphi+h\varphi^{p-1}=\lambda f\varphi^{\alpha-1}
\end{equation*}
on a finite graph $G$ with $\alpha\geq p>1$. This work is inspired by Grigor'yan, Lin and Yang's pioneer paper \cite{GLY,GLY'}, where they studied similar equations on finite or locally finite graphs.

\section{Settings and main results}
Let $G=(V,E)$ be a finite graph, where $V$ denotes the vertex set and $E$ denotes the edge set. Fix a vertex measure $\mu:V\rightarrow(0,+\infty)$ and an edge measure $\omega:E\rightarrow(0,+\infty)$ on $G$. The edge measure $\omega$ is assumed to be symmetric, that is, $\omega_{ij}=\omega_{ji}$ for each edge $i\thicksim j$.

Denote $C(V)$ as the set of all real functions defined on $V$, then $C(V)$ is a finite dimensional linear space with the usual function additions and scalar multiplications. For any $p>1$, the $p$-th discrete graph Laplacian $\Delta_p:C(V)\rightarrow C(V)$ is
\begin{equation*}
\Delta_pf_i=\frac{1}{\mu_i}\sum\limits_{j\thicksim i}\omega_{ij}|f_j-f_i|^{p-2}(f_j-f_i)
\end{equation*}
for any $f\in C(V)$ and $i\in V$. $\Delta_p$ is a nonlinear operator when $p\neq2$.

\label{sect-main-result}
\begin{theorem}\label{thm-main}
Let $G=(V,E)$ be a finite connected graph. Given $h, f\in C(V)$ with $f>0$. Assume $\alpha\geq p>1$. Then the following $p$-th Yamabe equation
\begin{equation}
\Delta_p\varphi+h\varphi^{p-1}=\lambda f\varphi^{\alpha-1}
\end{equation}\label{def-Yamabe-equ-p}
on $G$ always has a positive solution $\varphi$ for some constant $\lambda\in\mathds{R}$.
\end{theorem}

Taking $p=2$, we get the following
\begin{corollary}\label{croll-main}
Let $G=(V,E)$ be a finite connected graph. Given $h, f\in C(V)$ with $f>0$. Assume $\alpha>2$. Then the following Yamabe equation
\begin{equation}\label{def-Yamabe-equ-2}
\Delta\varphi+h\varphi=\lambda f\varphi^{\alpha-1}
\end{equation}
on $G$ always has a positive solution $\varphi$ for some constant $\lambda\in\mathds{R}$.
\end{corollary}

\begin{remark}
Grigor'yan, Lin and Yang \cite{GLY'} established similar results for the following equation
\begin{equation}\label{equ-gly-2}
-\Delta u+hu=|u|^{\alpha-2}u,\;\;\alpha>2
\end{equation}
on a finite graph under the assumption $h>0$. They show that the above equation (\ref{equ-gly-2}) always has a positive solution. They also studied the following equation
\begin{equation}\label{equ-gly-p}
-\Delta_p u+h|u|^{p-2}u=f(x,u),\;\;p>1
\end{equation}
and established some existence results under certain assumptions of $f(x,u)$. However, it is remarkable that their $\Delta_p$ considered in the equation (\ref{equ-gly-p}) is different with ours when $p\neq2$. It is also remarkable that our Theorem \ref{thm-main} doesn't require $h>0$.
\end{remark}

\section{Proofs of theorem \ref{thm-main}}
\label{sect-preliminary-lemma}
\subsection{Sobolev embedding}
For any $f\in C(V)$, define an integral of $f$ over $V$ with respect to the vertex weight $\mu$ by
$$\int_Vfd\mu=\sum\limits_{i\in V}\mu_if_i.$$
Set $\mathrm{Vol}(G)=\int_Vd\mu$. Similarly, for any function $g$ defined on the edge set $E$, we define an integral of $g$ over $E$ with respect to the edge weight $\omega$ by
$$\int_Egd\omega=\sum\limits_{i\thicksim j}\omega_{ij}g_{ij}.$$
Specially, for any $f\in C(V)$,
$$\int_E|\nabla f|^pd\omega=\sum\limits_{i\thicksim j}\omega_{ij}|f_j-f_i|^p,$$
where $|\nabla f|$ is defined on the edge set $E$, and $|\nabla f|_{ij}=|f_j-f_i|$ for each edge $i\thicksim j$. Next we consider the Sobolev space $W^{1,\,p}$ on the graph $G$. Define
$$W^{1,\,p}(G)=\left\{u\in C(V):\int_E|\nabla\varphi|^pd\omega+\int_V|\varphi|^pd\mu<+\infty\right\},$$
and
$$\|u\|_{W^{1,\,p}(G)}=\left(\int_E|\nabla\varphi|^pd\omega+\int_V|\varphi|^pd\mu\right)^{\frac{1}{p}}.$$
Since $G$ is a finite graph, then
$W^{1,\,p}(G)$ is exactly $C(V)$, a finite dimensional linear space. This implies the following Sobolev embedding:

\begin{lemma}\label{lem-Sobolev-embedding}(Sobolev embedding)
Let $G=(V,E)$ be a finite graph. The Sobolev space $W^{1,\,p}(G)$ is pre-compact. Namely, if $\{\varphi_n\}$ is bounded in $W^{1,\,p}(G)$, then there exists some $\varphi\in W^{1,\,p}(G)$ such that up to a subsequence, $\varphi_n\rightarrow\varphi$ in $W^{1,\,p}(G)$.
\end{lemma}
\begin{remark}
The convergence in $W^{1,\,p}(G)$ is in fact pointwise convergence.
\end{remark}
\subsection{Proofs step by step}
We follow the original approach pioneered by Yamabe \cite{Yamabe}. Denote an energy functional
\begin{equation}
I(\varphi)=\left(\int_E|\nabla \varphi|^pd\omega-\int_Vh\varphi^pd\mu\right)\left(\int_Vf\varphi^{\alpha} d\mu\right)^{-\frac{p}{\alpha}},
\end{equation}
where $\varphi\in W^{1,\,p}(G)$, $\varphi\geq 0$ and $\varphi\not\equiv0$. Define
\begin{equation}
\beta=\inf \big\{I(\varphi): \varphi\geq0,\;\varphi\not\equiv0\big\}.
\end{equation}
We shall find a solution to (\ref{def-Yamabe-equ-p}) step by step as follows.\\

\textbf{Step 1}. $I(\varphi)$ is bounded below for all $\varphi\geq0$, $\varphi\not\equiv0$. Hence $\beta\neq-\infty$ and $\beta\in\mathds{R}$. In fact, it's easy to see
$$0<\left(\int_Vf\varphi^{\alpha} d\mu\right)^{\frac{p}{\alpha}}\leq f_M^{\frac{p}{\alpha}}\left(\int_V\varphi^{\alpha} d\mu\right)^{\frac{p}{\alpha}}=f_M^{\frac{p}{\alpha}}\|\varphi\|_{\alpha}^p,$$
where $f_M=\max\limits_{i\in V}f_i>0$. Hence
\begin{equation}\label{equ-right}
\left(\int_Vf\varphi^{\alpha} d\mu\right)^{-\frac{p}{\alpha}}\geq f_M^{-\frac{p}{\alpha}}\|\varphi\|_{\alpha}^{-p}>0.
\end{equation}
Similarly, we also have
$$-\int_Vh\varphi^{p}d\mu\geq (-h)_m\int_V\varphi^{p}d\mu=(-h)_m\|\varphi\|_{p}^{p},$$
where $(-h)_m=\min\limits_{i\in V}(-h_i)$.
Then it follows
\begin{equation}\label{equ-left}
\int_E|\nabla \varphi|^pd\omega-\int_Vh\varphi^pd\mu\geq(-h)_m\|\varphi\|_{p}^{p}.
\end{equation}
By (\ref{equ-right}) and (\ref{equ-left}), we get
\begin{equation*}
I(\varphi)\geq(-h)_m\|\varphi\|_{p}^{p}f_M^{-\frac{p}{\alpha}}\|\varphi\|_{\alpha}^{-p},
\end{equation*}
and further
\begin{equation}\label{equ-I-fai-1}
I(\varphi)\geq\big((-h)_m\wedge 0\big)\|\varphi\|_{p}^{p}f_M^{-\frac{p}{\alpha}}\|\varphi\|_{\alpha}^{-p},
\end{equation}
where $(-h)_m\wedge 0$ is the minimum of $(-h)_m$ and $0$. Since $\alpha\geq p$, then
\begin{equation}\label{equ-half-1}
0<\|\varphi\|_{p}^{p}\leq\left(\int_V\left(\varphi^p\right)^{\frac{\alpha}{p}}d\mu\right)^{\frac{p}{\alpha}}
\left(\int_V1^{\frac{\alpha}{\alpha-p}}d\mu\right)^{\frac{\alpha-p}{\alpha}}
=\|\varphi\|_{\alpha}^{p}\mathrm{Vol}(G)^{1-\frac{p}{\alpha}},
\end{equation}
which leads to
\begin{equation}\label{equ-I-fai-2}
0<\|\varphi\|_{p}^{p}\|\varphi\|_{\alpha}^{-p}\leq\mathrm{Vol}(G)^{1-\frac{p}{\alpha}}.
\end{equation}
Thus by (\ref{equ-I-fai-1}) and (\ref{equ-I-fai-2}), we obtain
\begin{equation}\label{equ-I-fai-final}
I(\varphi)\geq\big((-h)_m\wedge 0\big)f_M^{-\frac{p}{\alpha}}\mathrm{Vol}(G)^{1-\frac{p}{\alpha}}=C_{\alpha,p,h,f,G},
\end{equation}
where $C_{\alpha,p,h,f,G}\leq0$ is a constant depending only on the information of $\alpha$, $p$, $h$, $f$ and $G$. Note that the information of $G$ contains $V$, $E$, $\mu$ and $\omega$. Hence $I(\varphi)$ is bounded below by a universal constant. \\

\textbf{Step 2}. There exists a $\hat{\varphi}\geq0$, such that $\beta=I(\hat{\varphi})$. To find such $\hat{\varphi}$, we choose $\varphi_n\geq0$, satisfying
$$\int_Vf\varphi_n^{\alpha}d\mu=1$$
and
$$I(\varphi_n)\rightarrow\beta$$
as $n\rightarrow\infty$. We may well suppose $I(\varphi_n)\leq 1+\beta$ for all $n$. Note
$$1=\int_Vf\varphi_n^{\alpha}d\mu\geq f_m\int_V\varphi_n^{\alpha}d\mu=f_m\|\varphi_n\|_{\alpha}^{\alpha},$$
where $f_m=\min\limits_{i\in V}f_i$. Hence
\begin{equation}\label{equ-half-2}
\|\varphi_n\|_{\alpha}^{p}\leq f_m^{-\frac{p}{\alpha}}.
\end{equation}
Denote $|h|_M=\max\limits_{i\in V}|h_i|$, then by (\ref{equ-half-1}) and (\ref{equ-half-2}), we obtain
\begin{equation*}
\begin{aligned}
\|\varphi_n\|^p_{W^{1,\,p}(G)}=&\int_E|\nabla\varphi|^pd\omega+\int_V|\varphi|^pd\mu\\
=&\;I(\varphi_n)+\int_Vh\varphi_n^{p}d\mu+\|\varphi_n\|_{p}^{p}\\
\leq&\;1+\beta+(1+|h|_M)\|\varphi_n\|_{p}^{p}\\
\leq&\;1+\beta+(1+|h|_M)\mathrm{Vol}(G)^{1-\frac{p}{\alpha}}\|\varphi_n\|_{\alpha}^{p}\\
\leq&\;1+\beta+(1+|h|_M)\mathrm{Vol}(G)^{1-\frac{p}{\alpha}}f_m^{-\frac{p}{\alpha}},
\end{aligned}
\end{equation*}
which implies that $\{\varphi_n\}$ is bounded in $W^{1,\,p}(G)$. Therefore by Lemma \ref{lem-Sobolev-embedding}, there exists some $\hat{\varphi}\in C(V)$ such that up to a subsequence, $\varphi_n\rightarrow \hat{\varphi}$
in $W^{1,\,p}(G)$. We may well denote this subsequence as $\varphi_n$. Note $\varphi_n\geq0$ and $\int_Vf\varphi_n^{\alpha}d\mu=1$, let $n\rightarrow+\infty$, we know $\hat{\varphi}\geq0$ and $\int_Vf\hat{\varphi}^{\alpha}d\mu=1$. This implies that $\hat{\varphi}\not\equiv0$. Since the energy functional $I(\varphi)$ is continuous, we have $\beta=I(\hat{\varphi})$.\\

\textbf{Step 3}. $\hat{\varphi}>0$.

Calculate the Euler-Lagrange equation of $I(\varphi)$, we get
\begin{equation}\label{equ-Euler-Lagrange}
\frac{d}{dt}\Big|_{t=0}I(\varphi+t\phi)=-p\left(\int_Vf\varphi^{\alpha} d\mu\right)^{-\frac{p}{\alpha}}\int_V\left(\Delta_p\varphi+h\varphi^{p-1}-\lambda_{\varphi}f\varphi^{\alpha-1}\right)\phi d\mu,
\end{equation}
where
\begin{equation}\label{def-lamda-fai}
\lambda_{\varphi}=-\frac{\int_E|\nabla \varphi|^pd\omega-\int_Vh\varphi^pd\mu}{\int_Vf\varphi^{\alpha} d\mu}
\end{equation}
for any $\varphi\geq0$, $\varphi\not\equiv0$. Thus
\begin{equation}\label{equ-gradient-I-i}
\frac{\partial I}{\partial \varphi_i}=-p\mu_i(\Delta_p\varphi_i+h\varphi^{p-1}_i-\lambda_{\varphi}f_i\varphi^{\alpha-1}_i)\left(\int_Vf\varphi^{\alpha} d\mu\right)^{-\frac{p}{\alpha}}.
\end{equation}
Note the graph $G$ is connected, if $\hat{\varphi}>0$ is not satisfied, since $\hat{\varphi}\geq0$ and not identically zero, then there is an edge $i\thicksim j$, such that $\hat{\varphi}_i=0$, but $\hat{\varphi}_j>0$. Now look at $\Delta_p\hat{\varphi}_i$,
$$\Delta_p\hat{\varphi}_i=\frac{1}{\mu_i}\sum\limits_{k\thicksim i}\omega_{ik}|\hat{\varphi}_k-\hat{\varphi}_i|^{p-2}(\hat{\varphi}_k-\hat{\varphi}_i)>0.$$
Therefore by (\ref{equ-gradient-I-i}), we have
$$\frac{\partial I}{\partial \varphi_i}\Big|_{\varphi=\hat{\varphi}}=-p\mu_i\Delta_p\hat{\varphi}_i\left(\int_Vf\hat{\varphi}^{\alpha} d\mu\right)^{-\frac{p}{\alpha}}<0.$$
Recall we had proved that $\hat{\varphi}$ is the minimum value of $I(\varphi)$, hence there should be
$$\frac{\partial I}{\partial \varphi_i}\Big|_{\varphi=\hat{\varphi}}\geq0,$$
which is a contradiction. Hence $\hat{\varphi}>0$.\\

\textbf{Step 4}. $\hat{\varphi}$ satisfied the equation (\ref{def-Yamabe-equ-p}), that is
\begin{equation}\label{equ-final}
\Delta_p\hat{\varphi}+h\hat{\varphi}^{p-1}=\lambda_{\hat{\varphi}} f\hat{\varphi}^{\alpha-1},
\end{equation}
where $\lambda_{\hat{\varphi}}$ is defined according to (\ref{def-lamda-fai}). Because $I(\varphi)$ attains its minimum value at $\hat{\varphi}$, which lies in the interior of $\{\varphi\in C(V):\varphi\geq0\}$, so
$$\frac{d}{dt}\Big|_{t=0}I(\hat{\varphi}+t\phi)=0$$
for all $\phi\in C(V)$. This leads to (\ref{equ-final}).\\

\noindent \textbf{Acknowledgements:} The author would like to thank Professor Gang Tian and Yanxun Chang for constant encouragement. The author would also like to thank Dr. Wenshuai Jiang, Xu Xu for many helpful conversations. The research is supported by National Natural Science Foundation of China under Grant No.11501027, and Fundamental Research Funds for the Central Universities (Nos. 2015JBM103, 2014RC028, 2016JBM071 and 2016JBZ012).

Huabin Ge: hbge@bjtu.edu.cn

Department of Mathematics, Beijing Jiaotong University, Beijing 100044, P.R. China
\end{document}